\documentclass[10pt]{article}
\usepackage[T1]{fontenc}
\usepackage[utf8]{inputenc}

\usepackage{lmodern}
\usepackage{xspace}

\usepackage{amsfonts,amsmath,amssymb, amsthm, mathtools}
\usepackage{dsfont} 
\usepackage{algorithmicx,algpseudocode,algorithm}

\usepackage[usenames,dvipsnames,table]{xcolor}
\usepackage[backref,colorlinks,citecolor=blue,bookmarks=true,linktocpage]{hyperref}
\usepackage[numbered]{bookmark}
\usepackage{fullpage}

\usepackage[shortlabels]{enumitem}
  \setitemize{noitemsep,topsep=3pt,parsep=2pt,partopsep=2pt} \setenumerate{itemsep=1pt,topsep=2pt,parsep=2pt,partopsep=2pt}
  \setdescription{itemsep=1pt}

\usepackage{mleftright}
\usepackage{aliascnt}
\makeatletter
\@ifundefined{theorem}{\theoremstyle{plain} \newtheorem{theorem}{Theorem}\newaliascnt{coro}{theorem}
  	  
  	\aliascntresetthe{coro}
  	\newaliascnt{lem}{theorem}
  		
  	\aliascntresetthe{lem}
  	\newaliascnt{clm}{theorem}
  		
	\aliascntresetthe{clm}
	\newaliascnt{fact}{theorem}
 	 	\newtheorem{fact}[theorem]{Fact}
	\aliascntresetthe{fact}
  	
  \newaliascnt{prop}{theorem}
  		\newtheorem{proposition}[prop]{Proposition}
	\aliascntresetthe{prop}
	\newaliascnt{conj}{theorem}
  		
	\aliascntresetthe{conj}
  \theoremstyle{remark}

  \theoremstyle{definition} \newaliascnt{defn}{theorem}
 		 
 	 \aliascntresetthe{defn}
}{}
\makeatother

\newenvironment{proofof}[1]{\begin{proof}[Proof of {#1}]}{\end{proof}}

     \newcommand{\eqdef}{\coloneqq}

 \newcommand{\distribs}[1]{\Delta\!\left(#1\right)}

\newcommand{\bigO}[1]{{O\mleft( #1 \mright)}}

\newcommand{\bigTheta}[1]{{\Theta\mleft( #1 \mright)}}
\newcommand{\bigOmega}[1]{{\Omega\mleft( #1 \mright)}}

 			\newcommand{\indicSet}[1]{\mathds{1}_{#1}}                                              \newcommand{\indic}[1]{\indicSet{\left\{#1\right\}}}                                             
\newcommand{\dtv}{\operatorname{d}_{\rm TV}}
\newcommand{\kl}{\operatorname{KL}}
\newcommand{\dhell}{\operatorname{d_{\rm{}H}}}
\newcommand{\hellinger}[2]{{\dhell\mleft({#1, #2}\mright)}}
\newcommand{\kldiv}[2]{{\kl\mleft({#1 \,\|\, #2}\mright)}}
\newcommand{\kolmogorov}[2]{{\operatorname{d_{\rm{}K}}\mleft({#1, #2}\mright)}}
\newcommand{\totalvardistrestr}[3][]{{\dtv^{#1}\mleft({#2, #3}\mright)}}
\newcommand{\totalvardist}[2]{\totalvardistrestr[]{#1}{#2}}
\newcommand{\chisquare}[2]{{\chi^2\mleft({#1 \mid\mid #2}\mright)}}

\newcommand\restr[2]{{\left.\kern-\nulldelimiterspace #1 \vphantom{\big|} \right|_{#2} }}

\newcommand{\proba}{\Pr}
\newcommand{\probaOf}[1]{\proba\!\left[\, #1\, \right]}

\newcommand{\expect}[1]{\mathbb{E}\!\left[#1\right]}

\newcommand{\shortexpect}{\mathbb{E}}
\newcommand{\var}{\operatorname{Var}}

\newcommand{\binomial}[2]{\ensuremath{\operatorname{Bin}\!\left( #1, #2 \right)}}

\newcommand{\norm}[1]{\lVert#1{\rVert}}
\newcommand{\normone}[1]{{\norm{#1}}_1}
\newcommand{\normtwo}[1]{{\norm{#1}}_2}
\newcommand{\norminf}[1]{{\norm{#1}}_\infty}
\newcommand{\abs}[1]{\left\lvert #1 \right\rvert}

\newcommand{\R}{\ensuremath{\mathbb{R}}\xspace}

\newcommand{\lp}[1][1]{\ell_{#1}}

\makeatletter
  \AtBeginDocument{
  \begingroup
  \toks@={}\toksdef\toks@@=2 \toks@@={}\long\def\@ReturnFiFi#1#2\fi\fi{\fi\fi#1}\def\scan@author#1#2 \and#3\@nil{\ifx\\#3\\\ifcase#1 \toks@={#2}\else
      \ifnum#1>1 \toks@=\expandafter{\the\expandafter\toks@\expandafter,\expandafter\space
          \the\toks@@
        }\fi
      \toks@=\expandafter{\the\toks@\space and #2}\fi
    \else
      \ifcase#1 \toks@={#2}\@ReturnFiFi{\scan@author1#3\@nil
        }\else
        \ifnum#1>1 \toks@=\expandafter{\the\expandafter\toks@\expandafter,\expandafter\space
            \the\toks@@
          }\fi
      \toks@@={#2}\@ReturnFiFi{\scan@author2#3\@nil
      }\fi
  \fi
  }\expandafter\expandafter\expandafter\scan@author
  \expandafter\expandafter\expandafter0\expandafter\@author\space\and\@nil
  \edef\x{\endgroup
  \noexpand\hypersetup{pdfauthor={\the\toks@}}}\x
  }
\makeatother
 
\newcommand{\dst}{\varepsilon}
\newcommand{\ab}{k}
\newcommand{\ns}{n}

\usepackage{filecontents}
\usepackage{utopia}

\title{A short note on learning discrete distributions}
\date{February, 2020}
\author{Cl\'{e}ment L. Canonne\thanks{The latest version of this note can be found at \href{https://github.com/ccanonne/probabilitydistributiontoolbox}{github.com/ccanonne/probabilitydistributiontoolbox}.}}

\begin{document}
\maketitle

\begin{abstract}
The goal of this short note is to provide simple proofs for the ``folklore facts'' on the sample complexity of learning a discrete probability distribution over a known domain of size $\ab$ to various distances $\dst$, with error probability $\delta$.
\end{abstract}

\tableofcontents\bigskip

\paragraph{Notation.}
For a given distance measure $\operatorname{d}$, we write $\Phi(\operatorname{d},\ab,\dst,\delta)$ for the sample complexity of learning discrete distributions over a known domain of size $\ab$, to accuracy $\dst>0$, with error probability $\delta\in(0,1]$. As usual, asymptotics will be taken with regard to $\ab$ going to infinity, $\dst$ going to $0$, and $\delta$ going to $0$, in that order.
 Without loss of generality, we hereafter assume the domain is the set $[\ab]\eqdef \{1,\dots,\ab\}$.

\section{Total variation distance}

Recall that $\totalvardist{p}{q} = \sup_{S\subseteq [\ab]} (p(S)-q(S)) = \frac{1}{2}\normone{p-q}\in[0,1]$ for any $p,q\in\distribs{[\ab]}$. 
\begin{theorem}\label{theo:learning:tv}
  $\Phi(\dtv,\ab,\dst,\delta) = \bigTheta{\frac{\ab+\log(1/\delta)}{\dst^2}}$.
\end{theorem}
\noindent We focus here on the upper bound. The lower bound can be proven, e.g., via Assouad's lemma (for the $\ab/\dst^2$ term), and from the hardness of estimating the bias of a coin ($\ab=2$) with high probability (for the $\log(1/\delta)/\dst^2$ term).
\begin{proof}[First proof]
Consider the empirical distribution $\tilde{p}$ obtained by drawing $\ns$ independent samples $s_1,\dots,s_\ns$ from the underlying distribution $p\in\distribs{[\ab]}$:
\begin{equation}\label{def:empirical}
\tilde{p}(i) = \frac{1}{\ns} \sum_{j=1}^\ns \indic{s_j=i}, \qquad i\in [\ab]
\end{equation}
\begin{itemize}
  \item First, we bound the \emph{expected} total variation distance between $\tilde{p}$ and $p$, by using $\lp[2]$ distance as a proxy:
\[
    \expect{ \totalvardist{p}{\tilde{p}} }
    =\frac{1}{2}\expect{ \normone{p-\tilde{p}}} 
    =\frac{1}{2}\sum_{i=1}^\ab\expect{ \abs{p(i)-\tilde{p}(i)}}
    \leq\frac{1}{2}\sum_{i=1}^\ab\sqrt{\expect{ (p(i)-\tilde{p}(i))^2} }
\]
the last inequality by Jensen. But since, for every $i\in[\ab]$, $\ns\tilde{p}(i)$ follows a $\binomial{\ns}{p(i)}$ distribution, we have
$\expect{ (p(i)-\tilde{p}(i))^2} = \frac{1}{\ns^2}\var[\ns\tilde{p}(i)] = \frac{1}{\ns}p(i)(1-p(i))$, from which
\[
    \expect{ \totalvardist{p}{\tilde{p}} } \leq\frac{1}{2\sqrt{\ns}}\sum_{i=1}^\ab\sqrt{p(i)} \leq \frac{1}{2}\sqrt{\frac{\ab}{\ns}}
\]
the last inequality this time by Cauchy--Schwarz. Therefore, for $\ns\geq \frac{\ab}{\dst^2}$ we have $\expect{ \totalvardist{p}{\tilde{p}} }\leq \frac{\dst}{2}$.

  \item Next, to convert this expected result to a \emph{high probability} guarantee, we apply McDiarmid's inequality to the random variable $f(s_1,\dots,s_\ns) \eqdef \totalvardist{p}{\tilde{p}}$, noting that changing any single sample cannot change its value by more than $c\eqdef 1/\ns$:
\[
    \probaOf{ \abs{f(s_1,\dots,s_\ns) - \expect{f(s_1,\dots,s_\ns)}} \geq \frac{\dst}{2} } \leq 2e^{-\frac{2\left(\frac{\dst}{2}\right)^2}{\ns c^2}} = 2e^{-\frac{1}{2}\ns\dst^2}
\]
and therefore as long as $\ns\geq \frac{2}{\dst^2}\ln\frac{2}{\delta}$, we have $\abs{f(s_1,\dots,s_\ns) - \expect{f(s_1,\dots,s_\ns)}} \leq \frac{\dst}{2}$ with probability at least $1-\delta$. 
\end{itemize}
Putting it all together, we obtain that $\totalvardist{p}{\tilde{p}} \leq \dst$ with probability at least $1-\delta$, as long as $\ns\geq \max\left( \frac{\ab}{\dst^2},\frac{2}{\dst^2}\ln\frac{2}{\delta} \right)$.
\end{proof}

\begin{proof}[Second proof -- the ``fun'' one]
Again, we will analyze the behavior of the empirical distribution $\tilde{p}$ over $\ns$ i.i.d. samples from the unknown $p$ (cf.~\eqref{def:empirical}) -- because it is simple, efficiently computable, and \emph{it works}.  Recalling the definition of total variation distance, note that $\totalvardist{p}{\tilde{p}} > \dst$ literally means there exists a subset $S\subseteq [\ab]$ such that $\tilde{p}(S) > p(S) + \dst$. There are $2^\ab$ such subsets, so\dots{} let us do a union bound.

Fix any $S\subseteq[\ab]$. We have
\[
\tilde{p}(S) = \tilde{p}(i) \operatorname*{=}^{\eqref{def:empirical}} \frac{1}{\ns} \sum_{i\in S} \sum_{j=1}^\ns \indic{s_j=i}
\]
and so, letting $X_j \eqdef \sum_{i\in S}\indic{s_j=i}$ for $j\in [\ns]$, we have
$
\tilde{p}(S) = \frac{1}{\ns}\sum_{j=1}^\ns X_j
$ where the $X_j$'s are i.i.d. Bernoulli random variable with parameter $p(S)$. Here comes the Chernoff bound (actually, Hoeffding, the \emph{other} Chernoff):
\[
    \probaOf{ \tilde{p}(S) > p(S) + \dst } = \probaOf{ \frac{1}{\ns}\sum_{j=1}^\ns X_j > \expect{\frac{1}{\ns}\sum_{j=1}^\ns X_j} + \dst } \leq e^{-2\dst^2 \ns}
\]
and therefore $\probaOf{ \tilde{p}(S) > p(S) + \dst } \leq \frac{\delta}{2^\ab}$ for any $\ns\geq \frac{\ab\ln 2+\log(1/\delta)}{2\dst^2}$. A union bound over these $2^\ab$ possible sets $S$ concludes the proof:
\[
    \probaOf{ \exists S\subseteq [\ab] \text{ s.t. }\tilde{p}(S) > p(S) + \dst } \leq 2^\ab\cdot \frac{\delta}{2^\ab} = \delta
\]
and we are done. \emph{Badda bing badda boom}, as someone\footnote{John Wright.} would say.
\end{proof}

\section{Hellinger distance}

Recall that $\hellinger{p}{q} = \frac{1}{\sqrt{2}}\sqrt{\sum_{i=1}^\ab (\sqrt{p(i)}-\sqrt{q(i)})^2} = \frac{1}{\sqrt{2}}\normtwo{\sqrt{p}-\sqrt{q}}\in[0,1]$ for any $p,q\in\distribs{[\ab]}$. The Hellinger distance has many nice properties: it is well-suited to manipulating product distributions, its square is subadditive, and is always within a quadratic factor of the total variation distance; see, e.g.,~\cite[Appendix~C.2]{Canonne:15}.

\begin{theorem}\label{theo:learning:hellinger}
  $\Phi(\dhell,\ab,\dst,\delta) = \bigTheta{\frac{\ab+\log(1/\delta)}{\dst^2}}$.
\end{theorem}

This theorem is ``highly non-trivial'' to establish, however; for the sake of exposition, we will show increasingly stronger bounds, starting with the easiest to establish.
\begin{proposition}[Easy bound]\label{theo:learning:hellinger:easy}
  $\Phi(\dhell,\ab,\dst,\delta) = \bigO{\frac{\ab+\log(1/\delta)}{\dst^4}}$, and $\Phi(\dhell,\ab,\dst,\delta) = \bigOmega{\frac{\ab+\log(1/\delta)}{\dst^2}}$.
\end{proposition}
\begin{proof}
    This is immediate from~\autoref{theo:learning:tv}, recalling that $\frac{1}{2}\dtv^2\leq \dhell^2\leq \dtv$.
\end{proof}
\begin{proposition}[More involved bound]\label{theo:learning:hellinger:intermediate}
  $\Phi(\dhell,\ab,\dst,\delta) = \bigO{\frac{\ab}{\dst^2}+\frac{\log(1/\delta)}{\dst^4}}$.
\end{proposition}
\begin{proof}
    As for total variation distance, we consider the empirical distribution $\widehat{p}$ (cf.~\eqref{def:empirical}) obtained by drawing $\ns$ independent samples $s_1,\dots,s_\ns$ from $p\in\distribs{[\ab]}$.
    
    \begin{itemize}
      \item First, we bound the \emph{expected} squared Hellinger distance between $\widehat{p}$ and $p$: using the simple fact that
      $\hellinger{p}{q}^2 = 1-\sum_{i=1}^\ab \sqrt{p(i)q(i)}$ for any $p,q\in\distribs{[\ab]}$,
      \[
          \expect{ \hellinger{p}{\widehat{p}}^2 } = 1-\sum_{i=1}^\ab \sqrt{p(i)}\cdot \expect{\sqrt{\widehat{p}(i)}}\,.
      \]
      Now we would like to handle the square root inside the expectation, and \emph{of course} Jensen's inequality is in the wrong direction. However, for every nonnegative r.v. $X$ with positive expectation, letting $Y\eqdef X/\expect{X}$, we have that
      \begin{align*}
          \expect{\sqrt{X}} 
          &= \sqrt{\expect{X}}\cdot\expect{\sqrt{Y}}
          = \sqrt{\expect{X}}\cdot\expect{\sqrt{1+(Y-\expect{Y})})} \\
          &\geq \sqrt{\expect{X}} \mleft( 1+ \frac{1}{2}\expect{Y-\expect{Y}} - \frac{1}{2} \expect{(Y-\expect{Y})^2} \mright)
          = \sqrt{\expect{X}}\mleft(1-\frac{\var X}{2\expect{X}^2}\mright)
      \end{align*}
      where we used the inequality $\sqrt{1+x} \geq 1+\frac{x}{2}-\frac{x^2}{2}$, which holds for $x\geq -1$.\footnote{And is inspired by the Tayor expansion $\sqrt{1+x} = 1+\frac{x}{2} - \frac{x^2}{8} +o(x^2)$: there is \emph{some} intuition for it.}{} Since, for every $i\in[\ab]$, $\ns\widehat{p}(i)$ follows a $\binomial{\ns}{p(i)}$ distribution, we get
      \[
          \expect{ \hellinger{p}{\widehat{p}}^2 } \leq 1-\frac{1}{\sqrt{\ns}}\sum_{i=1}^\ab \sqrt{p(i)}\cdot\sqrt{\ns p(i)} \mleft(1-\frac{\ns p(i)(1-p(i))}{2\ns^2 p(i)^2}\mright)
          \leq 1 - \sum_{i=1}^\ab p(i) \mleft(1-\frac{1}{2\ns p(i)}\mright) = \frac{\ab}{2\ns}\,.
      \]
      Therefore, for $\ns\geq \frac{\ab}{\dst^2}$, we have $\expect{ \hellinger{p}{\widehat{p}}^2 }\leq \frac{\dst^2}{2}$.
      \item Next, to convert this expected result to a high probability guarantee, we \emph{would like} to apply McDiarmid's inequality to the random variable $f(s_1,\dots,s_\ns) \eqdef \hellinger{p}{\widehat{p}}^2$ as in the (first) proof of~\autoref{theo:learning:tv}; unfortunately, changing a sample can change the value by up to $c \approx 1/\sqrt{\ns}$, and McDiarmid will yield
only a vacuous bound.\footnote{Try it: it's a real bummer.}{} Instead, we will use a stronger, more involved concentration inequality:
      \begin{theorem}[{\cite[Theorem~8.6]{BLM:13}}]\label{theo:stronger:mcdiarmid}
          Let $f\colon \mathcal{X}^\ns \to \R$ be a measurable function, and let $X_1,\dots,X_\ns$ be independent random variables taking values in $\mathcal{X}$. Define $Z\eqdef f(X_1,\dots,X_\ns)$. Assume that there exist measurable functions $c_i\colon \mathcal{X}^\ns \to [0,\infty)$ such that, for all $x,y\in\mathcal{X}^\ns$,
          \[
              f(y) - f(x) \leq \sum_{i=1}^\ns c_i(x) \indic{x_i \neq y_i}\,.
          \]
          Then, setting $v \eqdef \shortexpect \sum_{i=1}^\ns c_i(x)^2$ and $v_\infty \eqdef \sup_{x\in\mathcal{X}^\ns} \sum_{i=1}^\ns c_i(x)^2$, we have, for all $t>0$,
          \[
              \probaOf{ Z \geq \expect{Z} + t } \leq e^{-\frac{t^2}{2v}}\,\qquad \probaOf{ Z \leq \expect{Z} - t } \leq e^{-\frac{t^2}{2v_\infty}}\,.
          \]
      \end{theorem}
      For our $f$ above, we have, for two any different $x, y \in [\ab]^\ns$ , that
      \begin{align*}
        f(y)-f(x) 
        &= \frac{1}{\sqrt{\ns}} \sum_{i=1}^\ab \sqrt{p(i)} \mleft( \sqrt{\sum_{j=1}^\ns \indic{x_j=i}} - \sqrt{\sum_{j=1}^\ns \indic{y_j=i}} \mright) \\
        &= \frac{1}{\sqrt{\ns}} \sum_{i=1}^\ab \sqrt{p(i)} \frac{\sum_{j=1}^\ns (\indic{x_j=i}-\indic{y_j=i})}{ \sqrt{\sum_{j=1}^\ns \indic{x_j=i}} + \sqrt{\sum_{j=1}^\ns \indic{y_j=i}} } \\
        &\leq \frac{1}{\sqrt{\ns}} \sum_{i=1}^\ab \sqrt{p(i)} \frac{\sum_{j=1}^\ns \indic{x_j=i}\indic{y_j\neq x_j}}{ \sqrt{\sum_{j=1}^\ns \indic{x_j=i}} }
        = \sum_{j=1}^\ns \underbrace{\sqrt{ \frac{ p_{x_j} }{ \ns\sum_{\ell=1}^\ns \indic{x_\ell=x_j} } }}_{c_j(x)} \cdot \indic{x_j\neq y_j}\,.
      \end{align*}
      In view of~\autoref{theo:stronger:mcdiarmid}, we then must evaluate
      \[
          v \eqdef \sum_{j=1}^\ns \expect{c_j(X)^2} = \frac{1}{\ns} \sum_{j=1}^\ns  \sum_{i=1}^\ab p(i)^2 \cdot \expect{ \frac{1}{1+ \sum_{\ell\neq j} \indic{X_\ell=i} } }
      \]
      where that last expectation is over $(x_\ell)_{\ell\neq j}$ drawn from $p^{\otimes(\ns-1)}$. Since $\sum_{\ell\neq j} \indic{X_\ell=i}$ is Binomially distributed with parameters $\ns-1$ and $p(i)$, we can use the simple fact that, for $N\sim\binomial{r}{\rho}$,
      \[
          \expect{\frac{1}{N+1}} = \frac{1-(1-\rho)^{r+1}}{\rho(r+1)} \leq \frac{1}{\rho(r+1)}
      \]
      to conclude that $v \leq \frac{1}{\ns^2} \sum_{j=1}^\ns \sum_{i=1}^\ab p(i)  = \frac{1}{\ns}$. By~\autoref{theo:stronger:mcdiarmid}, we obtain
      \[
          \probaOf{ \abs{f(s_1,\dots,s_\ns)-\expect{f(s_1,\dots,s_\ns)}} \geq \frac{\dst^2}{2} } \leq e^{-\frac{1}{8}\ns\dst^4}
      \]
      and therefore, as long as $\ns \geq \frac{8}{\dst^4}\ln\frac{1}{\delta}$, we have $\abs{f(s_1,\dots,s_\ns)-\expect{f(s_1,\dots,s_\ns)}} \leq \frac{\dst^2}{2}$ with probability at least $1-\delta$.
    \end{itemize}
    Putting it all together, we obtain that $\hellinger{p}{\widehat{p}}^2 \leq \dst^2$ with probability at least $1-\delta$, as long as $\ns \geq \max\mleft( \frac{\ab}{\dst^2}, \frac{8}{\dst^4}\ln\frac{1}{\delta}\mright)$.
\end{proof}
\noindent We finally get to the final, optimal bound:
\begin{proofof}{\autoref{theo:learning:hellinger}}
We will rely on a recent~--~and quite involved~--~result due to Agrawal [Agr19], analyzing the concentration of the empirical distribution $\widehat{p}$ in terms of its Kullback--Leibler (KL) divergence with regard to the true $p$,
\[
    \kldiv{\widehat{p}}{p} = \sum_{i=1}^\ab \widehat{p}(i) \ln \frac{\widehat{p}(i)}{p(i)} \in[0,\infty]\,.
\]
Observing that $\hellinger{p}{q}^2 \leq \frac{1}{2} \kldiv{p}{q}$ for any distributions $p, q$, the aforementioned result is actually stronger than what we need:
\begin{theorem}[{\cite[Theorem~1.2]{Agrawal:19}}]
  Suppose $\ns \geq \frac{\ab-1}{\alpha}$. Then
  \[
      \probaOf{ \kldiv{\widehat{p}}{p} \geq \alpha } \leq e^{-\ns\alpha}\mleft(\frac{e\alpha\ns}{\ab-1}\mright)^{\ab-1}\,.
  \]
\end{theorem}
\noindent In view of the above relation between Hellinger and KL, we will apply this convergence result with $\alpha \eqdef 2\dst^2$, obtaining
  \[
      \probaOf{ \hellinger{\widehat{p}}{p} \geq \dst } \leq e^{-2\ns\dst^2 + (\ab-1)\ln\frac{2e\ns\dst^2}{\ab-1}}\,.
  \]
  \begin{fact}
      For $\ns \geq \frac{15}{2e}\frac{\ab}{\dst^2}$, we have $(\ab-1)\ln\frac{2e\ns\dst^2}{\ab-1}\leq \ns\dst^2$. 
  \end{fact}
  \begin{proof}
  The conclusion is equivalent to $2e\cdot\ln\frac{2e\ns\dst^2}{\ab-1}\leq \frac{2e\ns\dst^2}{\ab-1}$, and thus follows from the fact that $x\geq 2e \ln x$ for $x\geq 15$.
  \end{proof}
  \noindent This fact implies that, for $\ns \geq \frac{15\ab}{2\dst^2}$, $\probaOf{ \hellinger{\widehat{p}}{p} \geq \dst } \leq e^{-\ns\dst^2}$. Overall, we obtain that $\hellinger{p}{\widehat{p}} \leq \dst$ with probability at least $1-\delta$ as long as $\ns \geq \max\mleft( \frac{15\ab}{2e\dst^2}, \frac{1}{\dst^2}\ln\frac{1}{\delta}\mright)$, as desired.
\end{proofof}

\section{$\chi^2$ and Kullback–-Leibler divergences}
In view of the previous section, some remarks on Kullback--Leibler (KL) and chi-squared ($\chi^2$) divergences. Recall their definition, for $p,q\in\distribs{[\ab]}$,
\[
    \kldiv{p}{q} = \sum_{i=1}^\ab p(i) \ln \frac{p(i)}{q(i)}\,,\qquad \chisquare{p}{q} = \sum_{i=1}^\ab \frac{(p(i)-q(i))^2}{q(i)}
\]
both taking values in $[0,\infty]$; as well as the chain of (easily checked) inequalities
\[
    2\totalvardist{p}{q}^2 \leq \kldiv{p}{q} \leq \chisquare{p}{q}\,,
\]
where the first one is Pinsker's. Importantly, KL and $\chi^2$ divergences are unbounded and asymmetric, so
the order of p and q matters \emph{a lot}: for instance, it is easy to show that, without strong assumptions on the
unknown distribution $p\in\distribs{[\ab]}$, the empirical estimator $\widehat{p}$ cannot achieve $\kldiv{p}{\widehat{p}} < \infty$ (resp., $\chisquare{p}{\widehat{p}} < \infty$) with any finite number of samples.\footnote{You can verify this: intuitively, the issue boils down to having to non-trivially learn even the elements of the support of $p$ that have arbitrarily small probability.}{} So, that's uplifting. (On the other hand, \emph{other} estimators than the empirical one, e.g., add-constant estimators, do provide good learning guarantees for those distance measures: see for instance~\cite{Kamath:15}).\smallskip

We are going to focus here on getting $\kldiv{\widehat{p}}{p}$ and $\chisquare{\widehat{p}}{p}$ down to $\dst$. Of course, in view of the inequalities above, the latter is at least as hard as the former, and a lower bound on both follows from that on $\dtv$: $\bigOmega{(\ab+\log(1/\delta))/\dst^2}$. And, behold! The result of Agrawal~\cite{Agrawal:19} used in the proof of~\autoref{theo:learning:hellinger} does provide the optimal upper bound on learning in KL divergence~--~and it is achieved by the usual suspect, the empirical estimator:
\begin{theorem}\label{theo:learning:kl}
  $\Phi(\kl,\ab,\dst,\delta) = \bigTheta{\frac{\ab+\log(1/\delta)}{\dst}}$, where by $\kl$ we refer to minimizing $\kldiv{\widehat{p}}{p}$.
\end{theorem}
\noindent The optimal sample complexity of learning in $\chi^2$ as a function of $\ab,\dst,\delta$, however, remains open.

\section{Briefly: Kolmogorov, $\lp[\infty]$, and $\lp[2]$ distances}
To conclude, let us briefly discuss three other distance measures: Kolmogorov (a.k.a., ``$\lp[\infty]$ between cumulative distribution functions''), $\lp[\infty]$, and $\lp[2]$:
\[
    \kolmogorov{p}{q} = \max_{i\in[\ab]} \abs{\sum_{j=1}^i p(i) - \sum_{j=1}^i q(i)} 
\] and
\[
    \lp[2](p,q) = \normtwo{p-q} = \sqrt{\sum_{i=1}^\ab (p(i)-q(i))^2},\qquad
    \lp[\infty](p,q) = \norminf{p-q} = \max_{i\in[\ab]} \abs{p(i)-q(i)}\,.
\]
A few remarks first. The Kolmogorov distance is actually defined for any distribution on $\R$, not necessarily discrete; one can equivalently define it as
$
\kolmogorov{p}{q} =\sup_{i} ( \shortexpect_p[\indicSet{(-\infty,i]}] - \shortexpect_q[\indicSet{(-\infty,i]}] )$. This has a nice interpretation: recalling the definition of TV distance, both are of the form $\sup_{f \in\mathcal{C}} (\shortexpect_p[f] - \shortexpect_q[f])
$ where $\mathcal{C}$ is a class of measurable functions.\footnote{Such metrics on the space of probability distributions are called \emph{integral probability metrics}.}{} For TV distance, $\mathcal{C}$ is the class of indicators of all measurable subsets; for Kolmogorov, this is the (smaller) class of indicators of intervals of the form $(-\infty, a]$. (For Wasserstein/EMD distance, this will be the class of continuous, $1$-Lipschitz functions.)

Second, because of the above, and also monotonicity of $\lp[p]$ norms, Cauchy--Schwarz, the fact that $\lp[1](p,q) = 2\totalvardist{p}{q}$, and elementary manipulations, we have
\begin{alignat*}{2}
    \lp[2](p,q)&\leq 2\totalvardist{p}{q} &&\leq \sqrt{\ab}\lp[2](p,q) \\
    \lp[\infty](p,q)&\leq \lp[2](p,q) &&\leq \sqrt{\lp[\infty](p,q)}, \\
   \frac{1}{2}\lp[\infty](p,q)&\leq \kolmogorov{p}{q} &&\leq \totalvardist{p}{q}\,.
\end{alignat*}
That can be useful sometimes. Now, I will only briefly sketch the proof of the next theorem: the lower bounds follow from the simple case $\ab=2$ (estimating the bias of a biased coin), the upper bounds are achieved by the empirical estimator (again). Importantly, the result for Kolmogorov distance \emph{still applies to continuous, arbitrary distributions}.
\begin{theorem}\label{theo:learning:kolmogorov}
  $\Phi(\operatorname{d_{\rm{}K}},\ab,\dst,\delta),\Phi(\lp[\infty],\ab,\dst,\delta),\Phi(\lp[2],\ab,\dst,\delta) = \bigTheta{\frac{\log(1/\delta)}{\dst^2}}$, independent of $\ab$.
\end{theorem}
\begin{proof}[Sketch]
The proof for Kolmogorov distance is the most involved, and follows from a \emph{very} useful and non-elementary theorem due to Dvoretzky, Kiefer, and Wolfowitz from 1956~\cite{DKW:56} (with the optimal constant due to Massart, in 1990~\cite{Massart:90}):
\begin{theorem}[DKW Inequality]
  Let $\hat{p}$ denote the empirical distribution on $\ns$ i.i.d. samples from $p$ (an \emph{arbitrary} distribution on $\R$). Then, for every $\dst > 0$,
  \[
      \probaOf{ \kolmogorov{\hat{p}}{p} > \dst } \leq 2e^{-2\ns \dst^2 }\,.
  \]
  Note, again, that this holds even if $p$ is a \emph{continuous} (or arbitrary) distribution on an unbounded support.
\end{theorem}
The proof for $\lp[\infty]$ just follows the Kolmogorov upper bound and the aforementioned inequality $\lp[\infty](p,q) \leq 2\kolmogorov{p}{q}$ (which hinges on the fact that $p(i) = \sum_{j=1}^i p(i) - \sum_{j=1}^{i-1} p(i)$ and the triangle inequality). Finally, the proof for $\lp[2]$ is a nice exercise involving analyzing the expectation of the $\lp[2]^2$ distance achieved by the empirical estimator, and McDiarmid's inequality.
\end{proof}

\paragraph{Acknowledgments.}
 I would like to thank Gautam Kamath and John Wright for suggesting {``someone should write this up as a note,''} and to Jiantao Jiao for discussions about the Hellinger case.

\bibliographystyle{alpha}
\bibliography{references-learning}

\begin{thebibliography}{KOPS15}

\bibitem[Agr19]{Agrawal:19}
Rohit Agrawal.
\newblock Multinomial concentration in relative entropy at the ratio of
  alphabet and sample sizes.
\newblock {\em CoRR}, abs/1904.02291, 2019.

\bibitem[BLM13]{BLM:13}
St{\'e}phane Boucheron, G{\'a}bor Lugosi, and Pascal Massart.
\newblock {\em Concentration inequalities: A nonasymptotic theory of
  independence}.
\newblock Oxford University Press, 2013.

\bibitem[Can15]{Canonne:15}
Cl{\'{e}}ment~L. Canonne.
\newblock {A Survey on Distribution Testing: Your Data is Big. But is it Blue?}
\newblock {\em Electronic Colloquium on Computational Complexity {(ECCC)}},
  22:63, 2015.

\bibitem[DKW56]{DKW:56}
Aryeh Dvoretzky, Jack Kiefer, and Jacob Wolfowitz.
\newblock Asymptotic minimax character of the sample distribution function and
  of the classical multinomial estimator.
\newblock {\em Ann. Math. Statist.}, 27:642--669, 1956.

\bibitem[KOPS15]{Kamath:15}
Sudeep Kamath, Alon Orlitsky, Dheeraj Pichapati, and Ananda~Theertha Suresh.
\newblock On learning distributions from their samples.
\newblock In {\em Proceedings of The 28th Conference on Learning Theory},
  volume~40 of {\em Proceedings of Machine Learning Research}, pages
  1066--1100. PMLR, 2015.

\bibitem[Mas90]{Massart:90}
Pascal Massart.
\newblock The tight constant in the {D}voretzky-{K}iefer-{W}olfowitz
  inequality.
\newblock {\em Ann. Probab.}, 18(3):1269--1283, 1990.

\end{thebibliography}
\end{document}